\documentclass{amsart}
\usepackage{amssymb}
\usepackage{epsfig}

\author{Thomas Lam }
\thanks{The author was supported in part by a Clay Institute Liftoff Fellowship.}
\address{Department of Mathematics, M.I.T., Cambridge, MA 02139}

\email{thomasl (at) math.mit.edu}
\date{July~16, 2005}

\theoremstyle{plain}
\newtheorem{theorem}{Theorem}
\newtheorem{thm}[theorem]{Theorem}
\newtheorem{prop}[theorem]{Proposition}

\newtheorem{lemma}[theorem]{Lemma}

\theoremstyle{definition}
\newtheorem{definition}[theorem]{Definition}
\newtheorem{example}[theorem]{Example}

\theoremstyle{remark}
\newtheorem{remark}{Remark}

\newcommand{\brac}[1]{\left(#1\right)}

\def\uqsln{U_q(\widehat{ \mathfrak{
sl}}_n)}
\def\ll{\lambda}

\def\Q{\mathbb Q}

\def\N{\mathbb N}

\def\wt{\mathrm{wt}}
\def\End{\mathrm{End}}

\def\ag{{\mathfrak g}}
\def\ug{U_q(\ag)}

\newcommand{\ip}[1]{\langle #1 \rangle}

\def\g{\mathcal{G}}

\def\m{\Phi}
\def\Z{\mathbb{Z}}

\def\p{\mathcal{P}}
\def\d{\delta}
\def\ll{\lambda}
\def\ff{\mathcal{F}}
\def\f{\mathbf{F}}

\title{A combinatorial generalization of the Boson-Fermion
correspondence}

\begin{document}
\begin{abstract}
We attempt to explain the ubiquity of tableaux and of Pieri and
Cauchy formulae for combinatorially defined families of symmetric
functions.  We show that such formulae are to be expected from
symmetric functions arising from representations of Heisenberg
algebras.  The resulting framework that we describe is a
generalization of the classical Boson-Fermion correspondence, from
which Schur functions arise.  Our work can be used to understand
Hall-Littlewood polynomials, Macdonald polynomials and Lascoux,
Leclerc and Thibon's ribbon functions, together with other new
families of symmetric functions.
\end{abstract}
\maketitle
\section{Introduction}
\label{sec:intro}

The classical Boson-Fermion correspondence is an isomorphism between
two representations of the Heisenberg algebra $H$: the Bosonic Fock
space $K[H_{-}]$ and the Fermionic Fock space $\ff^{(0)}$. It
identifies the Schur functions $s_\ll(x_1,x_2,\ldots)$ as the images
of the basis of semi-infinite wedges $v_{i_1} \wedge v_{i_2} \wedge
\cdots$ under this isomorphism.  The Boson-Fermion correspondence is
an important basic result in mathematical physics; see for
example~\cite{KR}.

The aim of this article is to replace the classical Fermionic Fock
space in the Boson-Fermion correspondence by another representation
of the Heisenberg algebra, and to obtain other interesting families
of symmetric functions instead of the Schur functions.  The
symmetric functions that we obtain have a tableaux-like definition,
and satisfy both Pieri-like identities and a Cauchy-like identity,
which we now explain.

\medskip

Let $\{F_\lambda(x_1,x_2,\ldots) \in \Lambda_K : \lambda \in S\}$ be
a family of symmetric functions with coefficients in a field $K$
(usually $\Q$, $\Q(q)$ or $\Q(q,t)$), where $S$ is some indexing
set.  Many important families of symmetric functions have the
following trio of properties.
\begin{enumerate}
\item
They can be expressed as the generating functions for a set of
``tableaux'', which gives the monomial expansion of $F_\lambda$:
\[
F_\lambda(x_1,x_2,\ldots) = \sum_{T} s(T) x^{\wt(T)},
\]
where the sum is over tableaux $T$ with ``shape'' $\lambda$. The
composition $\wt(T)$ is the {\it weight} of $T$ and $s(T) \in K$ is
some additional parameter associated to $T$.

%Often the set $S$ has a poset structure $(S, <)$ where each pair
%$(\ll,\mu)$ such that $\ll < \mu$ has been given a weighting
%$\wt(\ll,\mu) \in K$. Then the tableaux are chains $T = (\ll^{(0)}
%< \ll^{(1)} < \cdots < \ll^{(r)})$ in $S$ with weighting $\wt_T =
%\prod_{i=1}^r \wt(\ll^{(i-1)},\ll^{(i)})$.

\item Together with a closely related dual family $\{G_\lambda(x_1,x_2,\ldots) : \lambda \in
S\}$ of symmetric functions, they satisfy a Cauchy identity:
\[
\sum_{\lambda \in S} F_\lambda(x_1,x_2,\ldots)
G_\lambda(y_1,y_2,\ldots) = \prod_{i,j = 1}^\infty \left(b_0 + b_1
x_iy_j + b_2 (x_iy_j)^2 + \cdots \right),
\]
where the coefficients $b_i \in K$.
\item
They satisfy a Pieri formula:
\[
\tilde{h}_k(x_1,x_2,\ldots) F_\lambda(x_1,x_2,\ldots) = \sum_{\mu
\rightharpoonup_k \lambda} b_{\lambda,\mu} F_\mu(x_1,x_2,\ldots),
\]
where $k \in \Z$ is a positive integer, $\{\tilde{h}_1,\tilde{h}_2,
\ldots \} \in \Lambda_K$ is a sequence of symmetric functions and
$b_{\lambda,\mu} \in K$ are coefficients for each pair $\lambda,\mu$
satisfying some condition $\mu \rightharpoonup_k \lambda$.
\end{enumerate}

In all such cases that the author is aware of, the definition of a
tableaux involves the condition $\mu \rightharpoonup_k \lambda$ in
the Pieri formula. The simplest case is when $K = \Q$ and $F_\lambda
= s_\lambda$, the family of Schur functions. The indexing set $S =
\p$ is the set of partitions. The tableaux are usual semi-standard
Young tableaux $T$; the statistic $s(T)$ is equal to 1 and $\wt(T)$
is the usual weight associated to $T$.  The dual family $\{G_\lambda
= s_\lambda\}$ is equal to the Schur functions again and in the
Cauchy formula, all the coefficients $b_i = 1$. In the Pieri
formula, $\tilde{h}_k = h_k$ are the homogeneous symmetric
functions. The condition $\mu \rightharpoonup_k \lambda$ is that
$\mu/\ll$ is a horizontal strip of size $k$ and all the coefficients
$b_{\ll,\mu} = 1$.  Recall in particular that a semi-standard Young
tableaux is just a chain of partitions forming a sequence of
horizontal strips.

\medskip

Understanding the ubiquity of these three properties in families of
symmetric functions was one of the main aims of our work. Our main
result is as follows.  Given a representation $V$ of a Heisenberg
algebra $H$ with a distinguished basis $\{v_s \mid s \in S\}$,
together with a highest vector $v_b$ in $V$, we define a family
$F^V_{s}(x_1,x_2,\ldots)$ (and a dual family $G^V_{s}$) of symmetric
functions which satisfy a generalized Boson-Fermion correspondence.
The definition of $F^V_{s}$ is tableaux-like: for example it gives
the monomial expansion of $F^V_{s}$.  We show in addition that
$F^V_s$ satisfy a Pieri rule and a Cauchy identity.  Examples of
symmetric functions that can be obtained in this way include the
Schur functions, Schur $Q$-functions, Hall-Littlewood functions and
Macdonald polynomials; see~\cite{Mac}.

The motivating example for us was actually a family
$\g_\ll(x_1,x_2,\ldots;q)$ of $q$-symmetric functions defined by
Lascoux, Leclerc and Thibon~\cite{LLT} combinatorially via ribbon
tableaux and algebraically using the action of the Heisenberg
algebra on the Fock space of the quantized affine algebra $\uqsln$.
In~\cite{Lam} we studied the $\g_\ll$ in analogy with Schur
functions and discovered ribbon Cauchy and Pieri identities.  The
current work is an attempt to understand this in a more systematic
and general framework.  As an application, we now give natural
generalizations of the functions $\g_\ll$ to Fock spaces of other
types and also to higher level Fock spaces.  By our main result,
these new symmetric functions satisfy Cauchy and Pieri rules as
well, and will be the subject of later work.

Our Pieri and Cauchy formulae depend heavily on a sequence $a_i$ of
parameters defining the relations of the Heisenberg algebra $H =
H[a_i]$ (see Section~\ref{sec:classicalBF}).  On the other hand, as
an abstract algebra, the Heisenberg algebra does not depend on the
$a_i$ (as long as they are non-zero).  Thus it is not clear
immediately which sequences $a_i$ would lead to an interesting
theory of symmetric functions.

\medskip
Our work is also closely related to more combinatorial work of
Fomin~\cite{Fom1,Fom2,Fom} and of Bergeron and Sottile~\cite{BS}.
Fomin is mostly concerned with Schensted correspondences and the
Cauchy identities while Bergeron and Sottile's work has led to
relations with non-commutative symmetric functions and to Hopf
algebras.  It seems that an interesting non-commutative version of
our theory also exists, though we have not attempted to make this
precise in the present article.

It would be most interesting to investigate other families of
symmetric functions which arise using our correspondence from other
representations of Heisenberg algebras which occur naturally.

\medskip

We now briefly describe the organization of the rest of the paper.
In Section~\ref{sec:schur}, we review the theory of Schur functions
and symmetric functions.  In Section~\ref{sec:classicalBF}, we
describe the classical Boson-Fermion correspondence.  In
Section~\ref{sec:symHei}, we explain how to obtain symmetric
functions from representations of Heisenberg algebras.  In
Section~\ref{sec:genBF}, we prove our generalized Boson-Fermion
correspondence.  In Section~\ref{sec:piericauchy}, we prove Pieri
and Cauchy formulae for our families of symmetric functions.  In
Section~\ref{sec:converse}, we prove a partial converse to the
theorems of Sections~\ref{sec:genBF} and~\ref{sec:piericauchy}.  In
Section~\ref{sec:BFexamples}, we give a series of examples beginning
with Schur functions, Macdonald polynomials and the behavior when
taking direct sums or tensor products of representations.  We then
explain the example of Lascoux, Leclerc and Thibon's ribbon
functions studied in~\cite{LLT,Lam}.  Finally, we explain how to
generalise ribbon functions to other types and higher levels,
following work of Kashiwara, Miwa, Petersen and Yung~\cite{KMPY} and
Takemura and Uglov~\cite{TU}.

\medskip
{\bf Acknowledgements.}  I thank Sergey Fomin for discussions
related to this work.  This work is part of my Ph.D Thesis at
M.I.T., written under the guidance of Richard Stanley.

\section{Schur functions}
\label{sec:schur} We will follow mostly the notation of~\cite{Mac}.
Let $K$ be a field with characteristic 0.  Let $\Lambda_K$ denote
the ring of symmetric functions over $K$.  The ring $\Lambda_K$
should be thought of as the ring of formal power series in countably
many variables $x_1,x_2,\ldots$, of bounded degree.  If the variable
set is important then we write $\Lambda_K(X)$ or $\Lambda_K(Y)$.  We
will let $h_1,h_2,\ldots$ denote the {\it homogeneous symmetric
functions} and $p_1,p_2,\ldots$ denote the {\it power sum symmetric
functions}. Each of these sets forms a set of algebraically
independent generators for $\Lambda_K$.

Let $\p$ denote the set of partitions.  Let $\ll = (\ll_1 \geq \ll_2
\geq \cdots \geq \ll_l > 0) \in \p$ be a partition.  The {\it size}
$|\ll|$ of $\ll$ is equal to $\ll_1 + \cdots + \ll_l$ and we write
$\ll \vdash |\ll|$.  We also write $l(\lambda) = l$.  We generally
do not distinguish between a partition $\ll$ and its Young diagram
$D(\ll)$.  If $D(\mu) \subset D(\lambda)$ then $\lambda/\mu$ is a
{\it skew shape} with size $|\lambda/\mu| = |\lambda| - |\mu|$.

We let $h_\ll := h_{\ll_1}h_{\ll_2}\cdots h_{\ll_l}$ and $p_\ll :=
p_{\ll_1}p_{\ll_2}\cdots p_{\ll_l}$.  The sets $\{h_\ll : \ll \in
\p\}$ and $\{p_\ll : \ll \in \p\}$ are bases of $\Lambda_K$. The
homogeneous symmetric functions and the power sum symmetric
functions are related by the formula
\begin{equation}\label{eq:ph}
h_n = \sum_{\ll \vdash n} z_\ll^{-1} p_\ll
\end{equation}
where $z_\ll = 1^{m_1(\ll)}m_1(\ll)! 2^{m_2(\ll)}m_2(\ll)! \cdots $
and $m_i(\ll)=|\{j \mid \ll_j = i\}|$.

The {\it monomial symmetric functions} are denoted $m_\ll$ and the
{\it Schur functions} are denoted $s_\ll$.  The Schur functions (and
more generally skew Schur functions) are the generating functions of
{\it Young tableaux}:
\begin{equation}
\label{eq:young} s_\ll(x_1,x_2,\ldots) = \sum_{T} x^{\wt(T)},
\end{equation}
where the sum is over all semistandard Young tableaux $T$ of shape
$\ll$.  Alternatively, $s_\ll = \sum_{\mu} K_{\ll\mu} m_\mu$ where
the {\it Kostka number} $K_{\ll\mu}$ is equal to the number of
semistandard Young tableaux of shape $\ll$ and weight $\mu$.  For
the purposes of this paper, a Young tableaux $T$ of shape $\ll$
should be thought of as a chain of partitions $T = (\emptyset =
\ll^0 \subset \ll^1 \subset \cdots \subset \ll^l = \ll)$ such that
each skew shape $\ll^i/\ll^{i-1}$ is a {\it horizontal strip}.  A
horizontal strip is a skew shape containing at most one box in each
column.  The weight of $T$ is then the composition $\wt(T) =
(|\ll^1/\ll^0|,|\ll^2/\ll^1|,\ldots,|\ll^l/\ll^{l-1}|)$.  Similarly
a Young tableaux of skew shape $\lambda/\mu$ is a chain of
partitions $(\mu = \ll^0 \subset \ll^1 \subset \cdots \subset \ll^l
= \ll)$.

The Schur functions satisfy the following Pieri formula, which
describes how to write the product of a Schur function and a
homogeneous symmetric function in terms of Schur functions:
\begin{equation}
\label{eq:pieri} h_k s_\ll = \sum_{\mu \rightharpoonup_k \ll} s_\mu,
\end{equation}
where here $\mu \rightharpoonup_k \ll$ means that the skew shape
$\mu/\ll$ is a horizontal strip of size $k$.

The Schur functions also satisfy the following Cauchy formula, which
holds within the ring $\Lambda_K(X) \widehat{\otimes}_K
\Lambda_K(Y)$, which is the completion of the tensor product of two
copies of the symmetric functions.
\begin{equation}
\label{eq:cauchy} \prod_{i,j}\frac{1}{1-x_iy_j} = \sum_{\ll}
s_\ll(x_1,x_2,\ldots)s_\ll(y_1,y_2,\ldots).
\end{equation}

The ring of symmetric functions $\Lambda_K$ possesses a bilinear
symmetric form $\ip{.,.}: \Lambda_K \times \Lambda_K \rightarrow K$
given by $\ip{s_\lambda,s_\mu} = \delta_{\lambda\mu}$, or
alternatively by $\ip{p_\lambda,p_\mu} =
\delta_{\lambda\mu}z_\lambda$. This inner product is known as the
{\it Hall inner product}.  If $f \in \Lambda_K$ then $f^\perp \in
\End(\Lambda_K)$ denotes the linear operator adjoint to
multiplication by $f$.  As a particular case $p_k^\perp = k
\frac{\partial}{\partial p_k}$ where the differential operator acts
on symmetric functions written as polynomials in the power sum
symmetric functions.

\section{The classical Boson-Fermion correspondence}
\label{sec:classicalBF} Let $K$ be a field with characteristic 0.
The Heisenberg algebra $H = H[a_i]$ denotes the associative algebra
over $K$ with 1 generated by $\{B_k: k \in \Z\backslash \{0\}\}$
satisfying
\[
[B_k,B_l] = l\cdot a_l \cdot \d_{k,-l},
\]
for some non-zero parameters $a_l \in K$ satisfying $a_l = -
a_{-l}$.  As an abstract algebra, $H$ does not depend on the choice
of the elements $a_l$, since the generators $B_k$ can be re-scaled
to force $a_l = 1$.  However, we shall be concerned with
representations of $H$, and some choices of the generators $B_k$
will be more natural.

Let $K[H_{-}] = K[B_{-1},B_{-2},\ldots]$ denote the Bosonic Fock
space representation of $H$.  The action of the Heisenberg algebra
on $K[B_{-1},B_{-2},\ldots]$ is determined by letting $B_k$ act by
multiplication for $k < 0$ and setting $B_k \cdot 1 = 0$ for $k >
0$.

One can identify $K[B_{-1},B_{-2},\ldots]$ with the algebra
$\Lambda_K$ of symmetric functions over $K$ by identifying $B_{-k}$
with $a_k p_k$ for $k > 0$.  The action of $H$ on $\Lambda_K$ is
then given by
\[
B_k \longmapsto \begin{cases}a_{-k} p_{-k} &\mbox{for $k \leq -1$,}\\
k \frac{\partial}{\partial p_k} &\mbox{for $k \geq 1$.}
\end{cases}
\]

We will need the following standard lemma later.
\begin{lemma}
\label{lem:bkcommute} Let $k \geq 1$ be an integer and $\ll$ be a
partition. Then
\[
B_{-k} B_{\ll} = {k} a_{k}m_{k}(\ll)B_{\mu} + B_{\ll} B_{-k},
\]
where $\mu$ is $\ll$ with one less part equal to $k$. If $m_{k}(\ll)
= 0$ then the first term is just 0.
\end{lemma}

If $V$ is a representation of $H$, then a vector $v \in V$ is
called a {\it highest weight vector} if $B_k \cdot v = 0$ for $k >
0$.  The following result is well known.  See for
example~\cite[Proposition 2.1]{KR}.

\begin{prop}
\label{prop:hei} Let $V$ be an irreducible representation of $H$
with non-zero highest weight vector $v \in V$.  Then there exists a
unique isomorphism of $H$-modules $\phi: V \tilde{\longrightarrow}
K[B_{-1},B_{-2},\ldots]$ such that $\phi(v) = 1$.
\end{prop}

For the remainder of this section we assume that $H = H[1]$ is given
by the parameters $a_l = 1$ for $l \geq 1$ and $a_l = -1$ for $l
\leq -1$. Let $W = \oplus_{j \in \Z} Kv_j$ be an
infinite-dimensional vector space with basis $\{v_j : j \in Z\}$.
Let $\ff^{(0)}$ denote the vector space with basis given by
semi-infinite monomials of the form $v_{i_0} \wedge v_{i_{-1}}
\wedge \cdots$ where the indices satisfy:
\begin{enumerate}
\item[(i)] $i_0 > i_{-1} > i_{-2} > \cdots$
\item[(ii)] $i_k = k$ for $k$ sufficiently small.
\end{enumerate}
We will call $\ff^{(0)}$ the {\it Fermionic Fock space}.

\begin{remark}
Usually $\ff^{(0)}$ is considered a subspace of a larger space $\ff
= \oplus_{m \in \Z} \ff^{(m)}$.  The spaces $\ff^{(m)}$ are defined
as for $\ff^{(0)}$ with the condition (ii) replaced by the condition
(ii$^{(m)}$): $i_k = k-m$ for $k$ sufficiently small.
\end{remark}

Define an action of $H$ on $\ff^{(0)}$ by
\begin{equation}
\label{eq:boson}
 B_k \cdot (v_{i_0} \wedge v_{i_{-1}} \wedge \cdots) = \sum_{j \leq
0} v_{i_0} \wedge v_{i_{-1}} \wedge \cdots \wedge v_{i_{j-1}} \wedge
v_{i_j - k} \wedge v_{i_{j+1}} \wedge \cdots.
\end{equation}
The monomials are to be reordered according to the usual exterior
algebra commutation rules so that $v_{i_0} \wedge \cdots \wedge
v_{i_j} \wedge v_{i_{j+1}} \wedge \cdots = - v_{i_0} \wedge \cdots
\wedge v_{i_{j+1}} \wedge v_{i_{j}} \wedge \cdots$.  Thus the sum on
the right hand side of (\ref{eq:boson}) is actually finite so the
action is well defined.  One can check that we indeed do obtain an
action of $H$.  It is also not hard to see that the representation
of $H$ on $\ff^{(0)}$ is irreducible.

%This defines an action of $H$ with parameters $a_l = 1$ for $l
%\geq 1$ and $a_l = -1$ for $l \leq 1$.

%The space $\ff^{(0)}$ can be identified with $\f$ earlier in the
%thesis by $(v_{i_0} \wedge v_{i_{-1}} \wedge \cdots) \leftrightarrow
%(i_0 , i_{-1}+1,\ldots)$.

The vector $\bar{v} = v_0 \wedge v_{-1} \wedge \cdots \in
\ff^{(0)}$ is a highest weight vector for this action of $H$.  By
Proposition \ref{prop:hei}, there exists an isomorphism $\sigma:
\ff^{(0)} \rightarrow \Lambda_K$ sending $\bar{v} \mapsto 1$. An
algebraic version of the {\it Boson-Fermion correspondence}
identifies the image of $v_{i_0} \wedge v_{i_{-1}} \wedge \cdots$
under the isomorphism $\sigma$.

\begin{thm}[{\cite[Lecture 6]{KR}}]
\label{thm:BF} Let $\ll_k = i_{-k} + k$.  Then $\sigma(v_{i_0}
\wedge v_{i_{-1}} \wedge \cdots) = s_\ll$.
\end{thm}

In~\cite{KR}, this is called the ``second'' part of the
boson-fermion correspondence.  It is important in the study of a
family of non-linear differential equations known as the {\it
Kadomtzev-Petviashvili (KP) Hierarchy}.  The ``first'' part
consists of identifying the image of certain {\it vertex
operators} under $\sigma$.  The relationship between vertex
operators and symmetric function theory have been studied
previously in~\cite{Jin,Jin1,Mac}.

Our aim will be to generalise Theorem~\ref{thm:BF} to
representations of Heisenberg algebras with arbitrary parameters
$a_i \in K$.  We will see that the symmetric functions that one
obtains in this manner will always have a tableaux-like definition
and satisfy Pieri and Cauchy identities.  In our approach, we have
ignored the vertex operators, but it would be interesting to see how
they are related to our results.

\section{Symmetric functions from representations of Heisenberg
algebras} \label{sec:symHei}

Let $H = H[a_i]$ be the Heisenberg algebra with parameters $a_i \in
K$. Define $B_\ll := B_{\ll_1} B_{\ll_2} \cdots B_{\ll_{l(\ll)}}$.
Let $D_k := \sum_{\ll \vdash k} z_\ll^{-1} B_\ll$ and $U_k :=
\sum_{\ll \vdash k} z_\ll^{-1} B_{-\ll}$ where $z_\ll$ is as defined
in Section~\ref{sec:schur}.  Thus $B_\ll$ and $D_k$ are related in
the same way as $p_\ll$ and $h_k$ (see~(\ref{eq:ph})).

Similarly define $B_{-\ll}:=B_{-\ll_1} B_{-\ll_2} \cdots
B_{-\ll_{l(\ll)}}$ and $U_k:= \sum_{\ll \vdash k} z_\ll^{-1}
B_{-\ll}$.  Also let $S_\ll \in H$ be given by $S_\ll := \sum_\mu
z_\mu^{-1}\chi_\mu^\ll B_{-\mu}$ where the coefficients
$\chi_\mu^\ll$ are the characters of the symmetric group given by
$s_\ll = \sum_\mu z_\mu^{-1}\chi_\mu^\ll p_\mu$.

Let $V$ be a representation of $H$ with distinguished basis $\{v_s
: s \in S\}$ for some indexing set $S$.  For simplicity we will
assume that both $V$ and $S$ are $\Z$-graded so that $v_s \in V$
are homogeneous elements and $\deg(v_s) = \deg(s)$, and that each
graded component of $V$ is finite-dimensional.  We will also
assume that the action of $H$ is graded in the sense that
$\deg(B_k) = -mk$ for some $m \in \Z\backslash\{0\}$. Define an
inner product $\ip{.,.}: V \times V \rightarrow K$ on $V$ by
requiring that $\{v_s \mid s \in S\}$ forms an orthonormal basis,
so that $\ip{v_s,v_s'} = \delta_{ss'}$.

Let $s, t \in S$.  Define the generating functions
\begin{equation}
\label{eq:def} F^{V}_{s/t}(x_1,x_2,\ldots) = F_{s/t}(x_1,x_2,\ldots)
:= \sum_{\alpha} x^\alpha \ip{U_{\alpha_l} U_{\alpha_{l-1}}\cdots
U_{\alpha_1} \cdot t, s},
\end{equation} where the sum is over all compositions $\alpha =
(\alpha_1,\alpha_2, \ldots, \alpha_l)$.  Similarly define
\[ G^{V}_{s/t}(x_1,x_2,\ldots) = G_{s/t}(x_1,x_2,\ldots) = \sum_{\alpha} x^\alpha
\ip{D_{\alpha_l} D_{\alpha_{l-1}}\cdots D_{\alpha_1} \cdot s,
t}.\] Note that $F_{s/t}$ and $G_{s/t}$ are homogeneous with
degree $\frac{\deg(s) - \deg(t)}{m}$.  So in particular if
$\frac{\deg(s) - \deg(t)}{m}$ is negative or non-integral then the
generating functions are 0.  For convenience we let $U_\alpha:=
U_{\alpha_l} U_{\alpha_{l-1}}\cdots U_{\alpha_1}$ and $D_\alpha:=
D_{\alpha_l} D_{\alpha_{l-1}}\cdots D_{\alpha_1}$.

The above definitions should be thought of as a tableaux-like
definition, as the following example explains.

\begin{example}[Schur functions]
Let $H[a_i] = H[1]$ and $V = \ff^{(0)}$.  Set $S = \p$ and $v_\ll :=
v_{i_0} \wedge v_{i_{-1}} \wedge \cdots$, where $\ll_k = i_{-k} +
k$.  Then we have
\[
U_k \cdot v_\ll = \sum_{\mu \rightharpoonup_k \ll} v_\mu,
\]
where the sum is over all horizontal strips $\mu/\ll$ of size $k$.
So the definition (\ref{eq:def}) of $F_{s/t}$ reduces to
(\ref{eq:young}) -- the combinatorial definition of skew Schur
functions in terms of Young tableaux.
\end{example}

The following Proposition is immediate from the definition, since
$U_k$ commutes with $U_l$ and $D_k$ commutes with $D_l$ for all
$k,l \in \N$.

\begin{prop}
The generating functions $F_{s/t}$ and $G_{s/t}$ are symmetric
functions.
\end{prop}

As before, let $K[H_{-}] \subset H$ denote the subalgebra of $H$
generated by $\{B_k \mid k < 0\}$ and similarly define $K[H_+]
\subset H$. The definitions of $F_{s/t}$ and $G_{s/t}$ can be
rephrased in terms of the {\it Heisenberg-Cauchy elements}
$\Omega(H_{-},X)$ and $\Omega(H_{+},X)$ which lie in the completed
tensor products $K[H_{-}] \widehat{\otimes} \Lambda_K(X)$ and
$K[H_{+}] \widehat{\otimes} \Lambda_K(X)$ respectively:
\[
\Omega(H_{-},X) := \sum_\ll U_\ll \otimes m_\ll = \sum_\ll
z_\ll^{-1} B_{-\ll} \otimes p_\ll = \sum_\ll S_\ll \otimes s_\ll.
\]
The last two equalities follow from the classical Cauchy identity.
Also define $\Omega(H_{+},X) \in K[H_{+}] \hat{\otimes}
\Lambda_K(X)$ by $\Omega(H_{+},X) = \sum_\ll D_\ll \otimes m_\ll$.

Thus for example, one has
\[
F_{s/t}(x_1,x_2,\ldots) = \ip{\Omega(H_{-},X) \cdot v_t, v_s}
\]
and
\[
G_{s/t}(x_1,x_2,\ldots) = \ip{\Omega(H_{+},X) \cdot v_s, v_t}.
\]
One has in particular
\begin{equation}
\label{eq:Gs} G_{s/t}(x_1,x_2,\ldots) = \sum_\ll z_\ll^{-1}p_\ll
\ip{B_{\ll} \cdot v_s,v_t}.
\end{equation}

Now let $b \in S$ be such that $v_b$ is a highest weight vector
for $H$.  We will write $F_{s}:= F_{s/b}$ and $G_{s}:=G_{s/b}$.
The element $\Omega(H_{-},X) \cdot v_b \in V \hat{\otimes}
\Lambda_K(X)$  depends only on the choice of $v_b$. The symmetric
functions $F_s$ are the coefficients of $\Omega(H_{-},X) \cdot
v_b$ when it is written in the basis $\{v_s \mid s \in S\}$:
\[
\Omega(H_{-},X) \cdot v_b = \sum_s  v_s \otimes F_s(x_1,x_2,\ldots).
\]

\section{Generalization of Boson-Fermion correspondence}
\label{sec:genBF}  Let us suppose that $b \in S$ has been picked so
that $v_b \in V$ is a highest weight vector for $H$.  By
Proposition~\ref{prop:hei}, there is a canonical map of $H$-modules
$\phi: H \cdot b \rightarrow \Lambda_K$ sending $v_b \mapsto 1$. Our
choice of inner product for $V$ allows us to give a map $\m: V
\rightarrow \Lambda_K$.

\begin{theorem}[Generalized Boson-Fermion correspondence]
\label{thm:BFplus} The map $\m: V \rightarrow \Lambda_K$ given by
$v_s \mapsto G_s(x_1,x_2,\ldots)$ is a map of $H$-modules.
\end{theorem}

Recall that $B_{-k}$ acts on $\Lambda_K$ by multiplication by
$a_kp_k$ and $B_k$ acts as $k \frac{\partial}{\partial p_k}$, for
$k \geq 1$.
\begin{proof}
Let us calculate $B_l \cdot G_s$ and compare with $\m(B_l \cdot
v_s)$. Suppose first that $l < 0$ and let $k = -l$. Let $\ll$ be a
partition and let $\mu$ be $\ll$ with one less part equal to $k$.
If $\ll$ has no part equal to $k$, then $\mu$ can be any partition
in the following formulae.  First write $\ip{B_{\ll} B_l \cdot
v_s,v_b} = k a_k m_k(\ll) \ip{B_{\mu} v_s,v_b}$, using a slight
variation of Lemma~\ref{lem:bkcommute} for our $H$. Alternatively,
one can also compute
\begin{align*}
B_{\ll} B_l\cdot v_s &= B_\ll \sum_{c} \ip{B_l\cdot v_s,v_c} v_c =
\sum_{c,d} \ip{B_l\cdot v_s,v_c}\ip{B_\ll\cdot v_c,v_d}v_d
\end{align*}
so that taking the coefficient of $v_b$ we obtain \begin{equation}
\label{eq:use} k a_k m_k(\ll) \ip{B_{\mu} \cdot v_s,v_b} = \sum_c
\ip{B_l\cdot v_s,v_c}\ip{B_\ll\cdot v_c,v_b}.\end{equation}

Now,
\begin{align*} B_l \cdot G_s &= a_k p_k G_s
\\ &= a_k \sum_{\mu} z_\mu^{-1}  p_k p_\mu \ip{B_{\mu} \cdot v_s,
v_b} &\mbox{using (\ref{eq:Gs})},\\ &= \sum_{\ll} z_\ll^{-1} p_\ll
\left( \sum_c
\ip{B_l\cdot v_s,v_c}\ip{B_\ll\cdot v_c,v_b} \right) & \mbox{using (\ref{eq:use})}\\
& = \sum_c \ip{B_l\cdot v_s,v_c} \left( \sum_\ll z_\ll^{-1}
\ip{B_{\ll}\cdot v_c,v_b} \right) \\
&= \sum_c \ip{B_l \cdot v_s,v_c} G_c.
\end{align*}
This shows that $\m(B_l \cdot v_s) = B_l \cdot \m(v_s)$ for $l <
0$.

Now suppose $k > 0$, and let $\ll$ and $\mu$ be related as before.
Then
\begin{align*}
B_k \cdot G_s &= k \sum_{\ll} z_\ll^{-1} \frac{\partial}{\partial
p_k} p_\ll \ip{B_{\ll} \cdot v_s, v_b} \\
&=k \sum_{\ll} z_\ll^{-1} m_k(\ll) p_\mu \ip{B_\mu B_k \cdot v_s,v_b} \\
&=\sum_{\mu} z_\mu^{-1}p_\mu \ip{B_\mu \cdot \sum_{c} \ip{B_k\cdot
v_s,v_c}v_c, v_b} \\
&=\sum_c \ip{B_k \cdot v_s,v_c} \left( \sum_\mu z_\mu^{-1} p_\mu
\ip{B_\mu \cdot v_c, v_b} \right) \\
&= \sum_c \ip{B_k \cdot v_s,v_c} G_c.
\end{align*}

This completes the proof.
\end{proof}
When $V$ is irreducible, the map $\Phi$ does not depend on the
choice of basis, but does depend on $v_b$.  Since the degree
$\deg(v_b)$ part of $V$ is one dimensional, the image of $v \in V$
is given by the coefficient of the degree $\deg(v_b)$ part of
$\Omega(H_{+},X) \cdot v$.

If $V$ is not irreducible then the map depends on the inner product
$\ip{.,.}$ of $V$ (or equivalently, the choice of orthonormal
basis).

Note that a different action of $H$ on $\Lambda_K$ will allow us to
replace the family $G_s$ in Theorem~\ref{thm:BFplus} by $F_s$. More
precisely, one can define the adjoint action $\vartheta: H
\rightarrow \End(V)$ of $H$ on $V$ by letting the generators $B_k$
act according to the formula $ \ip{ \vartheta(B_k)\cdot v_{s'},v_s}
= \ip{v_{s'},B_{-k} \cdot v_s}$. With this new representation of $H$
on $V$, the roles of $G_s$ and $F_s$ are reversed.

\section{Pieri and Cauchy identities} \label{sec:piericauchy} Let
$h_k[a_i]$ denote the image $\theta(h_k)$ of $h_k$ under the algebra
homomorphism $\theta: \Lambda \rightarrow \Lambda_K$ given by
$\theta(p_k)= a_k p_k$. Also let $h_k\langle a_i \rangle$ denote the
image $\kappa(h_k)$ of $h_k$ under the map $\kappa:\Lambda_K
\rightarrow K$ given by $\kappa(p_k) = a_k$.   Note that if all
$\{a_i \mid i \geq 1\}$ are positive (rational) numbers then by
(\ref{eq:ph}) so are the numbers $h_k \langle a_i \rangle$.  Let
$h_k^\perp$ be the linear operator on $\Lambda_K$ which is adjoint
to multiplication by $h_k$ with respect to the Hall inner product.

\begin{thm}[Generalized Pieri Rule]
\label{thm:genPieri} Let $k \geq 1$.  The following identities
hold in $\Lambda_K$:
\[
h_k[a_i] G_s = \sum_{t} \ip{U_{k} \cdot s, t} G_t
\]
and
\[
h_k[a_i] F_s = \sum_{t} \ip{D_{k} \cdot t, s} F_t.
\]
The dual identities are:
\[
h_k^\perp G_s = \sum_{t} \ip{D_{k} \cdot s, t} G_t
\]
and
\[
h_k^\perp F_s = \sum_{t} \ip{U_{k} \cdot t, s} F_t.
\]
\end{thm}
\begin{proof}
Follows immediately from the definitions of $U_k,D_k$ and
$h_k[a_i]$ together with Theorem~\ref{thm:BFplus} and the comments
immediately after it.
\end{proof}

\begin{lemma}
\label{lem:du}   The following identity holds as elements of
$H[a_i]$:
\begin{equation}
\label{eq:du} D_b U_a = \sum_{j=0}^m h_j\langle a_i \rangle
U_{a-j}D_{b-j},
\end{equation}
where $m = \min(a,b)$.
\end{lemma}
\begin{proof}
By definition we need to show that $$ \brac{\sum_{\ll \vdash
b}z_\ll^{-1}B_\ll}\brac{ \sum_{\ll \vdash a} z_\ll^{-1}B_{-\ll}}  =
\sum_{j=0}^m h_j\langle a_i \rangle \brac{\sum_{\ll \vdash
 a-j}z_\ll^{-1} B_{-\ll}}\brac{ \sum_{\ll \vdash b-j}
 z_\ll^{-1}B_{\ll}}.
$$
Let $\mu$ and $\nu$ be partitions such that $|\mu| = a-j$ and $|\nu|
= b-j$. By (\ref{eq:ph}), the coefficient of $B_{-\mu} B_{\nu}$ on
the right hand side is equal to $z_\nu^{-1}z_\mu^{-1}\sum_{\ll
\vdash j} z_\ll^{-1}\theta(p_\ll)$. Let $\rho = \ll \cup \mu$ and
$\pi = \ll \cup \nu$.  We claim that the summand
$z_\nu^{-1}z_\mu^{-1}z_\ll^{-1}\theta(p_\ll)$ is the coefficient of
$B_{-\mu} B_\nu$ when applying $[B_k,B_l] = k a_k \delta_{k,-l}$
repeatedly to $z_\pi^{-1}z_\rho^{-1}B_{\pi}B_{-\rho}$. This is a
straightforward computation, counting the number of ways of picking
parts from $\rho$ and $\pi$ to make the partition $\ll$.

%Thus (\ref{eq:com}) implies Corollary \ref{cor:hhperp}, and since
%both the homogeneous and power sum symmetric functions generate the
%algebra of symmetric functions, Corollary~\ref{cor:hhperp} must be
%equivalent to (\ref{eq:com}).
\end{proof}

In fact the relation (\ref{eq:du}), together with the relations
$[U_k,U_l] = [D_k,D_l] = 0$ is equivalent to the defining relations
of the Heisenberg algebra $H[a_i]$.  This is because the sets $\{B_k
\mid k \neq 0\}$ and $\{U_k \mid k \geq 1\} \cup \{D_k \mid k \geq
1\}$ are both generators of $H[a_i]$.

%Define a map $\kappa: \Lambda \rightarrow K$ by $p_k \mapsto a_k$.

\begin{thm}[Generalized Cauchy Identity]
\label{thm:genCauchy} We have the following identity in the
completion of $\Lambda_K(X) \otimes \Lambda_K(Y)$:
\[
\sum_{s} F_s(x_1,x_2,\ldots) G_s(y_1,y_2,\ldots) =
\prod_{j,k}\left(1 + h_1\langle a_i\rangle x_jy_k +h_2\langle
a_i\rangle (x_jy_k)^2 + \cdots \right).
\]
More generally, let $r,t \in S$.  Then we have
\begin{multline*}
%\label{eq:genCauchy}
\sum_{s} F_{s/t}(x_1,x_2,\ldots) G_{s/r}(y_1,y_2,\ldots) = \\
\prod_{j,k}\left(1 +  h_1 \langle a_i\rangle x_jy_k +h_2\langle a_i
\rangle (x_jy_k)^2 + \cdots \right) \sum_s
F_{r/s}(x_1,x_2,\ldots)G_{t/s}(y_1,y_2,\ldots).
\end{multline*}
\end{thm}
\begin{proof}
Let $U(x):= 1+ \sum_{i > 0} U_i x^i$ and similarly $D(x):= 1+
\sum_{i > 0} D_i x^i$.  The identity of Lemma~\ref{lem:du} is
equivalent to
\[
D(y)U(x) = U(x)D(y)\left( 1 + h_1\langle a_i \rangle xy +h_2 \langle
a_i \rangle (xy)^2 + \cdots \right).
\]
Now notice that by definition we have $F_{s/t} = \ip{\cdots
U(x_3)U(x_2)U(x_1) \cdot v_t, v_s}$ and $G_{s/t} = \ip{\cdots
D(x_3)D(x_2)D(x_1) \cdot v_s, v_t}$.  The infinite products make
sense since in most factors we are picking the term equal to 1.
Thus
\begin{align*}
&\sum_{s} F_{s/t}(x_1,x_2,\ldots) G_{s/r}(y_1,y_2,\ldots)\\ &=
\ip{\cdots D(y_3)D(y_2)D(y_1) \cdots
U(x_3)U(x_2) U(x_1) \cdot v_t,v_r} \\
&= \prod_{i,j \geq 1}^\infty \left( 1 + h_1\langle a_i \rangle
x_iy_j +h_2\langle a_i \rangle (x_iy_j)^2 + \cdots \right) \\
& \hspace{30pt} \ip{\cdots U(x_3)U(x_2)
U(x_1) \cdots D(y_3)D(y_2)D(y_1) \cdot v_t,v_r} \\
&=\prod_{i,j \geq 1}^\infty \left( 1 + h_1\langle a_i \rangle x_iy_j
+h_2\langle a_i \rangle (x_iy_j)^2 + \cdots \right) \sum_s
G_{t/s}(y_1,y_2,\ldots)F_{r/s}(x_1,x_2,\ldots).
\end{align*}
These manipulations of infinite generating functions make sense
since they are well defined when restricted to a finite subset of
the variables $\{x_1,x_2,\ldots,y_1,y_2,\ldots\}$.
\end{proof}

\begin{remark}
It is not clear at this moment which sequences $a_i$ and which
representations of $H[a_i]$ would lead to interesting families of
symmetric functions. However, the following may be possible
indications:
\begin{itemize}
\item
Some kind of positivity for the coefficients $h_i\langle a_i
\rangle$; for example if $K = \Q(q)$ then we may want $h_i \langle
a_i \rangle$ to have positive coefficients when expanded as a power
series in $q$.
\item
A Pieri formula with very few non-zero or with positive
coefficients.  For example, we may want the coefficients $\ip{U_{k}
\cdot s, t}$ and $\ip{D_{k} \cdot t, s}$ to be positive in some
sense. This would imply that the definitions of $F_{t/s}$ and
$G_{t/s}$ would also have a positive monomial expansion.
\end{itemize}
\end{remark}

The results of this Section are related to results of
Fomin~\cite{Fom1,Fom2,Fom} and of Bergeron and Sottile~\cite{BS}.
Fomin studies combinatorial operators on posets and recovers Cauchy
style identities similar to ours.  His approach is more
combinatorial and he focuses on generalizing Schensted style
algorithms to these more general situations.  Bergeron and Sottile
have also made definitions similar to our $F_{s/t}$. Their interests
have been towards aspects related to Hopf algebras and
non-commutative symmetric functions; see also~\cite{BMSW,Ehr}.

\begin{remark} An interesting non-commutative version of our theory may exist,
where the Heisenberg algebra is replaced with an algebra $A =
\langle B_k \mid k \in \Z-\{0\} \rangle$ with relations
\begin{align*}
[B_k,B_l] &= 0 & \mbox{if $k$ and $l$ have opposite sign and $k \neq
-l$},\\ [B_{-k},B_{k}] &= ka_k.
\end{align*}
In this case, the generating functions $F_{s/t}$ and $G_{s/t}$ will
not be symmetric functions but instead be quasi-symmetric functions.
\end{remark}
\section{A partial converse}
\label{sec:converse} A partial converse to
Theorems~\ref{thm:genPieri} and~\ref{thm:genCauchy} exists.  In
other words, if a family of symmetric functions satisfies enough
properties, then one can conclude that they arise from a generalized
Boson-Fermion correspondence as in Theorem~\ref{thm:BFplus}.

Let $V$ be a $K$-vector space with a distinguished basis $\{v_s : s
\in S\}$.  In this section, suppose that $\{B_k' \in \End(V): k \in
\Z\backslash\{0\}\}$ are linear operators acting on $V$.  Suppose
further that $B_k$ and $B_l$ commute if $k$ and $l$ have the same
sign.  Let$D_k' := \sum_{\ll \vdash k} z_\ll^{-1} B_\ll'$ and $U_k'
:= \sum_{\ll \vdash k} z_\ll^{-1} B_{-\ll}'$.  Now we can define
$F'_{s/t}(x_1,x_2,\ldots) := \sum_{\alpha} x^\alpha
\ip{U'_{\alpha_l} U'_{\alpha_{l-1}}\cdots U'_{\alpha_1} \cdot t, s}$
and similarly for $G'_{s/t}$.

\begin{thm}
\label{thm:converse}  Let $\{a_k \in K \mid k \neq 0\}$ be a
sequence of non-zero parameters satisfying $a_k =a_{-k}$ and suppose
that $\{G'_s \mid s \in S\}$ are linearly independent. Then the
following are equivalent:
\begin{enumerate}
\item \label{it:BF} The operators $\{B'_k\}$ generate an action of the
Heisenberg algebra $H[a_i]$ with parameters $a_i$. \item
\label{it:pieri} The family $\{G'_s\}$ satisfies the conclusions of
Theorem~\ref{thm:genPieri}.
\item \label{it:cauchy} The families $\{G'_{s/t}\}$ and $\{F'_{s/t}\}$ satisfy the
conclusions of Theorem~\ref{thm:genCauchy}.
\end{enumerate}
\end{thm}

\begin{proof}
That (\ref{it:BF}) implies (\ref{it:pieri}) and (\ref{it:cauchy}) is
Theorems~\ref{thm:genPieri} and~\ref{thm:genCauchy}.

Now suppose (\ref{it:pieri}) holds.  Since the family $\{G'_s\}$ is
linearly independent, the action of $\{U'_k,D'_k\}$ on $V$ is
isomorphic to the action of $\{h_k[a_i], h_k^\perp\}$ on ${\rm
span}_K\{G'_s\}$ under the isomorphism $v_s \mapsto G'_s$. Thus the
action of the operators $B_k'$ on $V$ is isomorphic to the action of
$\{\theta(p_k), p_k^\perp\}$ on ${\rm span}_K\{G'_s\}$ and so
generate an action of $H[a_i]$.  Thus (\ref{it:pieri}) $\Rightarrow$
(\ref{it:BF}).

Now suppose (\ref{it:cauchy}) holds.  Then by the argument in the
proof of Theorem~\ref{thm:genCauchy}, we must have
\[
\ip{\left(D'(y)U'(x) - U'(x)D'(y)\left( 1 + h_1\langle a_i \rangle
xy +h_2\langle a_i \rangle(xy)^2 + \cdots \right) \right) \cdot v_t,
v_r} = 0
\]
for every $t, r \in S$.  This implies that \[D'(y)U'(x) =
U'(x)D'(y)\left( 1 + h_1\langle a_i \rangle xy +h_2\langle a_i
\rangle (xy)^2 + \cdots \right)\] so that we have \[D_b' U_a' =
\sum_{j=0}^m h_j\langle a_i \rangle U_{a-j}'D_{b-j}'.\] Now
reversing the argument in the proof of Lemma~\ref{lem:du}, we deduce
that $[B'_k,B'_{l}] = k a_k \delta_{k,-l}$.  So (\ref{it:cauchy})
$\Rightarrow$ (\ref{it:BF}).

\end{proof}

\section{Examples}
\label{sec:BFexamples}
\subsection{Schur functions}
If $K = \Q$ and $V = \ff^{(0)}$ and $H = H_{\text{Schur}} = H[1]$
acts as in Section~\ref{sec:classicalBF}, then
Theorem~\ref{thm:BFplus} is just Theorem~\ref{thm:BF}, where the
indexing set $S$ can be identified with the set of partitions $\p$.
In this case, the operators $B_k$ and $B_{-k}$ are adjoint with
respect to $\ip{.,.}$ and so $F_\ll = G_\ll = s_\ll$ for every
$\ll$.  The definition of $s_{\ll/\mu} = F_{\ll/\mu}$ in terms of
the operators $U_k$ is exactly the usual combinatorial definition of
skew Schur functions in terms of semistandard Young tableaux.  The
symmetric function $h_k[a_i] = h_k$ is the usual homogeneous
symmetric function and the coefficients $\ip{U_{k} \cdot \ll, \mu}$
are equal to 1 if $\mu/\ll$ is a horizontal strip of size $k$ and
equal to 0 otherwise.  The coefficients $h_i\langle a_i \rangle$ are
all equal to 1 and Theorem~\ref{thm:genCauchy} reduces to the usual
Cauchy identity.

\subsection{Direct sums}
\label{sec:direct} Let $V_1$ and $V_2$ be two representations of $H$
with distinguished bases $\{v_{s_1} : s_1 \in S_1\}$ and $\{v_{s_2}
: s_2 \in S_2\}$ respectively.  Then $V = V_1 \oplus V_2$ is a
representation of $H[a_i]$ with distinguished basis $\{v_{s} \mid s
\in S_1 \amalg S_2\}$.  If $s,t \in S_i$ for some $i$ then
$F^V_{s/t} = F^{V_i}_{s/t}$ otherwise if for example $s  \in S_1$
and $t \in S_2$ we have $F^V_{s/t} = 0$.  Thus the family of
symmetric functions that we obtain from $H[a_i]$ acting on $V$ is
the union of the families of symmetric functions we obtain from
$V_1$ and $V_2$.

\subsection{Tensor products}
\label{sec:tensor}

Let $V_1$ and $V_2$ be two representations of $H[a_i]$ with
distinguished bases $\{v_{s_1} : s_1 \in S_1\}$ and $\{v_{s_2} : s_2
\in S_2\}$ respectively, as before.  Then $V_1 \otimes V_2$ has a
distinguished basis $\{v_{s_1} \otimes v_{s_2} \mid s_1 \in S_1 \;
\text{and} \; s_2 \in S_2\}$.  Let the Heisenberg algebra
$\tilde{H}:= H[b_i]$ with generators $\tilde{B_k}$, where $b_i =
2a_i$, act on $V_1 \otimes V_2$ by defining the action of
$\tilde{B}_k$ by
\[
\tilde{B}_k \cdot v_1 \otimes v_2 = (B_k \cdot v_1) \otimes v_2 +
v_1 \cdot (B_k \cdot v_2).
\]
This action is natural when one views $\Lambda_K$ as a Hopf algebra.
The action of $\tilde{U}_k = \sum_{\lambda \vdash k} z_\lambda^{-1}
\tilde{B}_{-\ll}$ is given by
\[
\tilde{U}_k \cdot v_1 \otimes v_2 = \sum_{i = 0}^{k} (U_i \cdot
v_1) \otimes (U_{k-i} \cdot v_2)
\]
and similarly for $\tilde{D}_k$.  By definition, one sees that
$F^{V_1 \otimes V_2}_{s_1 \otimes s_2/ t_1 \otimes t_2} =
F^{V_1}_{s_1/t_1} F^{V_2}_{s_2/t_2}$ and similarly for the
$G$-functions. Thus the family of symmetric functions we obtain from
$V = V_1 \otimes V_2$ are pairwise products of the symmetric
functions we obtain from $V_1$ and $V_2$.

More generally, the tensor products $V_1 \otimes \cdots \otimes V_n$
lead to generating functions which are products $F^{V_1}_{s_1/t_1}
\cdots F^{V_n}_{s_n/t_n}$ of $n$ original generating functions.  We
will denote the Heisenberg algebra acting on this tensor product by
$H^{(n)}:=H[a_i^{(n)}]$.  The parameters are given by $a_i^{(n)} =
na_i$.

\subsection{Macdonald polynomials}
Let $K = \Q(q,t)$ and let $P_\ll(x_1,x_2,\ldots;q,t)$ and
$Q_\ll(x_1,x_2,\ldots;q,t)$ be the Macdonald polynomials introduced
in~\cite{Mac}.  Let $\ll = (\ll_1,\ll_2,\ldots)$ be a partition and
$s =(i,j) \in \ll$ be a square.  Then the arm-length of $s$ is given
by $a_\ll(s) = \ll_i - j$ and the leg-length of $s$ is given by
$l_\ll(s) = \ll'_j - i$. Now let $s$ be any square. Define
(\cite[Chapter VI, (6.20)]{Mac})
\[
b_\ll(s) = b_\ll(s;q,t) = \begin{cases}
\frac{1-q^{a_\ll(s)}t^{l_\ll(s) + 1}}{1-q^{a_\ll(s) +
1}t^{l_\ll(s)}} & \mbox{if $s \in \ll$,} \\
1 & \mbox{otherwise.} \end{cases}
\]
Now let $\ll/\mu$ be a horizontal strip.  Let $C_{\ll/\mu}$
(respectively $R_{\ll/\mu}$) denote the union of columns
(respectively rows) that intersect $\ll-\mu$.  Define
(\cite[Chapter VI, (6.24)]{Mac})
\[
\phi_{\ll/\mu} = \prod_{s \in C_{\ll/\mu}}
\frac{b_\ll(s)}{b_\mu(s)},
\]
and
\[
\psi_{\ll/\mu} = \prod_{s \in R_{\ll/\mu}-C_{\ll/\mu}}
\frac{b_\mu(s)}{b_\ll(s)}.
\]

Let $V_{\text{Mac}}$ denote the vector space over $K$ with
distinguished basis labeled by partitions.  Define operators $\{U_k,
D_k : k \in \Z_{>0}\}$ by:
\[
U_{k} \cdot \ll =  \sum_\mu \phi_{\mu/\ll} \mu, \;\;\;\;\; D_k \cdot
\ll = \sum_\mu \psi_{\ll/\mu}\mu,
\]
where the sums are over horizontal strips of size $|k|$. Then
$Q_{\ll/\mu} = F_{\ll/\mu}$ and $P_{\ll/\mu} = G_{\ll/\mu}$, so in
particular the operators $\{U_k \mid k \in \Z_{>0}\}$ commute and
so do the operators $\{D_k \mid k \in \Z_{>0}\}$.  Now we have
(\cite[Ex.7.6]{Mac})
\[
\sum_\rho Q_{\rho/\ll}(X;q,t) P_{\rho/\mu}(Y;q,t) = (\sum_\sigma
Q_{\mu/\sigma}(X;q,t)P_{\ll/\sigma}(Y;q,t)) \prod_{i,j}
\prod_{r=0}^\infty \frac{1-tx_iy_jq^r}{1-x_iy_jq^r}.
\]
The product $\prod_{r=0}^\infty \frac{1-ytq^r}{1-yq^r}$ can be
written as $\sum_{n \geq 0} g_n(1,0,0,\ldots;q,t)y^n$ where $g_n$
is given by (\cite[Chapter VI, (2.9)]{Mac})
\[
g_n(x_1,x_2,\ldots;q,t) = \sum_{\ll \vdash n} z_\ll(q,t)^{-1}
p_\ll(x_1,x_2,\ldots),
\]
where $z_\ll(q,t) = z_\ll \prod_{i=1}^{l(\ll)}
\frac{1-q^{\ll_i}}{1-t^{\ll_i}}$. Using Theorem~\ref{thm:converse},
we see that the operators $\{U_k,D_k \mid k \in \Z_{>0}\}$ generate
a copy of a Heisenberg algebra $H_{\text{Mac}}$.  A short
calculation shows that the parameters $a_k \in \Q(q,t)$ of this
Heisenberg algebra are given by $a_k = \frac{1-t^k}{1-q^k}$.  The
parameters $h_k \langle a_i \rangle$ are given by $h_k\langle a_i
\rangle = g_n(1,0,0,\ldots;q,t) = n \sum_{\ll \vdash n}
z_\ll(q,t)^{-1}$.

In fact Theorem~\ref{thm:converse} shows that the Pieri (and dual
Pieri) rule for Macdonald polynomials is equivalent to the
(generalized) Cauchy identity for Macdonald polynomials.

\begin{remark}
To obtain the Hall-Littlewood functions, one can just specialize
$q=0$ in the set up of this section.  However, to obtain the Schur
$P$ and $Q$-functions the further specialization $t = -1$ actually
causes some of the $a_i$ to be zero.  In this case, one should
actually consider the subalgebra of the Heisenberg algebra generated
by the generators $B_k$ where $k$ is odd.
\end{remark}

\subsection{Ribbon functions}  Let $n \geq 1$ be a positive integer and $K = \Q(q)$.
In~\cite{LLT}, a family of symmetric functions
$\{\g^{(n)}_\ll(x_1,x_2,\ldots;q)\}$ defined in terms of {\it ribbon
tableaux}, called {\it ribbon functions} or {\it LLT-polynomials},
were introduced.  These symmetric functions arise as the polynomials
$\{F^\f_s(x_1,x_2,\ldots)\}$ for the action of a Heisenberg algebra
$H[a_i]$ on the level one Fock space $\f$ of $\uqsln$.  This Fock
space $\f$ has a basis $|\lambda\rangle$ naturally labeled by
partitions.  The parameters are given by $a_i = \frac{1 -
q^{2nk}}{1-q^{2k}}$ and the action of $H[a_i]$, commuting with the
action of $\uqsln$, was discovered in~\cite{KMS}. The actions of the
generators $B_{-k}$ and $B_k$ of this Heisenberg algebra are adjoint
with respect to the inner product $\ip{|\lambda\rangle,|\mu\rangle}
= \delta_{\lambda\mu}$, and so the symmetric functions $F_\ll$ and
$G_\ll$ for this representation of $H[a_i]$ coincide. In~\cite{Lam},
a ribbon Cauchy and Pieri formula for the functions
$\g^{(n)}_\ll(X;q)$ was deduced from the action of $H[a_i]$ and this
is a special (in fact, motivating) case for
Theorems~\ref{thm:genPieri} and~\ref{thm:genCauchy}.

At $q = 1$, the Fock space $\f$ for $\uqsln$ should be thought of as
a sum of tensor products:
\begin{equation}
\label{eq:decomp} \f \cong \bigoplus_{\text{$n$-cores}}
(\ff^{(0)})^{\otimes n}
\end{equation}
where $\ff^{(0)}$ is the classical Fermionic Fock space described in
Section~\ref{sec:classicalBF}. Combinatorially, the decomposition
(\ref{eq:decomp}) is given by writing a partition in terms of its
$n$-core and its $n$-quotient; see~\cite{Mac}.  As shown in
subsection~\ref{sec:tensor}, the $F$-functions we obtain in this way
are products of $n$ of the $F$-functions for  $\ff^{(0)}$, that is,
(skew) Schur functions. This is simply the formula
$\g_\ll(x_1,x_2,\ldots;1) = s_{\ll^{(0)}}s_{\ll^{(1)}}\cdots
s_{\ll^{(n-1)}}$ observed in~\cite{LLT}.  In fact, the $q=1$
specialization corresponds to action of the Heisenberg algebra
commuting with the action of $\widehat{\mathfrak{sl}}_n$ on $\f$.

It would be interesting to see whether ribbon functions and
Macdonald polynomials can be combined by finding a deformation of
the action of $(H_{\text{Mac}})^{(n)}$ on $V_{\text{Mac}}^{\otimes
n}$.

\subsection{Ribbon functions for other types and other levels}
Theorem~\ref{thm:BFplus} allows us to define analogues of LLT's
ribbon functions $\g^{(n)}(x_1,x_2,\ldots;q)$ for other (quantized)
Fock spaces.

Kashiwara, Miwa, Petersen and Yung~\cite{KMPY} have defined (level
one) $q$-deformed Fock spaces for the affine algebras $A_n^{(1)}$,
$A_{2n}^{(2)}$, $B_n^{(2)}$, $A_{2n-1}^{(2)}$, $D_n^{(1)}$ and
$D_{n+1}^{(2)}$, using a sophisticated construction involving
perfect crystals.  Let $\Phi$ denote one of these root systems and
let $\ug$ be the corresponding quantum affine algebra.  Let
$\ff^\Phi$ be the corresponding $q$-deformed Fock space
of~\cite{KMPY}, which is defined over $K = \Q(q)$.  The space
$\ff^\Phi$ is equipped with an action of an Heisenberg algebra
$H[a_i^\Phi]$ commuting with the action of $\ug$, where the
parameters $a_i^\Phi$ are calculated in~\cite{KMPY}. The Fock space
$\ff^\Phi$ also has a standard basis indexed by certain
semi-infinite products of elements from a perfect crystal for $\ug$.
We will call this indexing set $S^\Phi$.  There is a distinguished
highest weight vector $v_b \in \ff^\Phi$ for some ``bottom element''
$b \in S^\Phi$.

\begin{definition}
Let $s \in S^\Phi$.  The {\it ribbon function of type $\Phi$} is
given by $\g^\Phi_s = F^{\ff^\Phi}_{s/b} \in \Lambda_K$.
\end{definition}

When $\Phi = A_{n-1}^{(1)}$, we recover LLT's ribbon functions
$\g^\Phi = \g^{(n)}(x_1,x_2,\ldots;q)$.  The functions
$\g^{(n)}(x_1,x_2,\ldots;q)$ have been found to be not only
interesting combinatorially (see~\cite{Lam, LLT}) but also to be
related to the global basis of the Fock space and to Kazhdan-Lusztig
polynomials (see~\cite{LT}). One should expect the symmetric
functions $\g^\Phi_s$ to be interesting as well.  Some work in this
direction can be found in~\cite{LamThesis} and will appear
separately.  Note that it is not known (but in some cases a
conjecture) that the action of the generators $B_k$ and $B_{-k}$ of
the Heisenberg algebra on $\ff^\Phi$ are adjoint.  This would imply
that $F^{\ff^\Phi}_{s/b} = G^{\ff^\Phi}_{s/b}$.
\medskip

In another direction, Takemura and Uglov~\cite{TU} have studied Fock
spaces $\f^{n,m}$ for the quantum affine algebra $\uqsln$ of level
$m$. These Fock spaces also possess a standard basis indexed by
partitions and an action of a Heisenberg algebra $H^{n,m}$ commuting
with the action of $\uqsln$.

\begin{definition}
Let $\ll \in \p$.  The {\it ribbon function of rank $n$ and level
$m$} is given by $\g^{(n,m)}_\ll = F^{\f^{n,m}}_{\ll/\emptyset} \in
\Lambda_K$.
\end{definition}

We have placed the parameters $n$ and $m$ together in the notation
since as explained in~\cite{TU} there is a level-rank duality in
this Fock space.  The case $m = 1$ reduces to LLT's ribbon
functions: $\g^{(n,1)}_\ll = \g^{(n)}_\ll$.  One should expect the
functions $\g^{(n,m)}_\ll$ to be interesting as well.  The
parameters $a_i$ for $H^{n,m}$ appear to have not yet been
calculated, though there are precise conjectures for their values.

\end{document}